# НЕРАВЕНСТВА ТИПА БЕРНШТЕЙНА ДЛЯ НЕПЕРИОДИЧЕСКИХ СПЛАЙНОВ В ПРОСТРАНСТВЕ $L_2$

В.Ф. Бабенко*'**, С.А. Спектор*.

Украина. *Днепропетровский национальный университет

** Институт прикладной математики и механики НАН Украины

Пусть $x \in R^1$, $\nu \in Z$ и $m = 0,1,2\ldots$ Функция $s(x) = \sum_\gamma c_\gamma N_m(x+\gamma)$, где

$$N_m(x) = \frac{1}{m!} \sum_{k=0}^{m+1} (-1)^k \binom{m+1}{k} (x-k)_+^m,$$

называется сплайном порядка $m$ минимального дефекта с узлами $\ell h$, $h > 0$, $\ell \in Z$, если:

1) $s(x)$ является полиномом с действительными коэффициентами степени $\leq n$ на каждом промежутке $(h(\ell-1), h\ell)$, $\ell \in Z$;

2) $s(x) \in C^{m-1}(R^1)$.

Совокупность всех таких сплайнов обозначим $S_{m,h}$. Через $\tilde{S}_{m,n}$ обозначим множество $2\pi$-периодических сплайнов порядка $m$, минимального дефекта, с узлами в точках $\frac{k\pi}{n}$, $n \in N$, $k = 0, \pm 1, \pm 2, \ldots$.

Во многих вопросах теории аппроксимации большую роль играют неравенства типа Бернштейна для периодических и непериодических сплайнов. Обзор известных точных неравенств типа Бернштейна для сплайнов из $\tilde{S}_{m,n}$ можно найти, например, в работе [1, 211-245].

Так, для периодических сплайнов неравенства вида

$$\left\|s_m^{(k)}\right\|_p \leq C_{m,p} n^k \left\|s_m\right\|_p$$

известны при $p = 1, 2, \infty$, а именно

$$\frac{\left\|s^{(k)}\right\|_1}{\left\|\varphi_{n,m-k}\right\|_1} \leq \frac{\|s\|_1}{\left\|\varphi_{n,m}\right\|_1}$$

$$\frac{\left\|s^{(k)}\right\|_2}{\left\|\varphi_{n,m-k}\right\|_2} \leq \frac{\|s\|_2}{\left\|\varphi_{n,m}\right\|_2}$$

$$\frac{\left\|s^{(k)}\right\|_{\infty}}{\left\|\varphi_{n,m-k}\right\|_{\infty}} \leq \frac{\|s\|_{C}}{\left\|\varphi_{n,m}\right\|_{C}},$$

где $n,m \in N$, $k = 1,...,m$, а $\varphi_{n,m}$ - стандартный совершенный сплайн. [1, 220, 229, 234]

Нами доказано следующее неравенство дающее оценку $k$-ой производной непериодического сплайна в пространстве $L_2$.

**Теорема.** Для всех функций $s \in S_{m,h}$, $\forall k \in Z$, $k \leq m$, $\forall h \in N$, справедливо

$$\left\|s^{(k)}\right\|_2 \leq (\pi h)^k \sqrt{\frac{K_{2(m-k)+1}}{K_{2m+1}}} \|s\|_2, \qquad (1)$$

где $K_m = \frac{4}{\pi} \sum_{\ell=0}^{\infty} \frac{(-1)^{\ell(m+1)}}{(1+2\ell)^{m+1}}$, $m = 0,1,2...$, – константы Фавара. [2, С. 64-65]

Константа в правой части неравенства (1) неулучшаема.

**Доказательство**:

Неравенство (1) достаточно доказать при $k=1$ и $\forall m$, а затем воспользоваться методом индукции.

Рассмотрим $s'(x) = \sum_{\gamma} c_{\gamma} N'_m(x+\gamma)$

Известно [3, 146], что $N'_m(x) = N_{m-1}(x) - N_{m-1}(x-1)$

Таким образом,

$$\|s'(x)\|_2^2 = \frac{1}{2\pi} \int_R \left| \sum_{\gamma} c_{\gamma} e^{-i\gamma x} \hat{N}_{m-1}(\omega)(1-e^{-i\omega}) \right|^2 d\omega =$$

$$\frac{1}{2\pi} \int_0^{2\pi} \sum_{\ell} |m_s(\omega)|^2 \left| \hat{N}_{m-1}(\omega+2\pi\ell)(1-e^{-i\omega}) \right|^2 d\omega =$$

$$\frac{1}{2\pi} \int_0^{2\pi} \frac{\sum_{\ell} |i\omega+2\pi\ell|^2 \left|\hat{N}_m(\omega+2\pi\ell)\right|^2}{\sum_{\ell}\left|\hat{N}_m(\omega+2\pi\ell)\right|^2} |m_s(\omega)|^2 \sum_{\ell}\left|\hat{N}_{m-1}(\omega+2\pi\ell)\right|^2 \left|1-e^{-i\omega}\right|^2 d\omega \leq$$

$$\max_{\omega} \frac{\left|\hat{N}_{m-1}(\omega+2\pi\ell)\right|^2 \left|1-e^{-i\omega}\right|^2}{\sum_{\ell}\left|\hat{N}_m(\omega+2\pi\ell)\right|^2} \|s(\omega)\|_2^2 \qquad (2)$$

Найдем

$$\max_{\omega} L(\omega) = \max_{\omega} \frac{\left|\hat{N}_{m-1}(\omega+2\pi\ell)\right|^2 \left|1-e^{-i\omega}\right|^2}{\sum_{\ell}\left|\hat{N}_m(\omega+2\pi\ell)\right|^2}.$$

Так как $\sum_{\ell}\left|\hat{N}_m(\omega+2\pi\ell)\right|^2 = \frac{1}{(2m-1)!}\prod_{\ell=1}^{m-1}\frac{1-2\lambda_{\ell}\cos\omega+\lambda_{\ell}^2}{|\lambda_{\ell}|}$

и $\sum_{\ell}\left|\hat{N}_{m-1}(\omega+2\pi\ell)\right|^2 = \frac{1}{(2m-3)!}\prod_{\ell=1}^{m-2}\frac{1-2\lambda'_{\ell}\cos\omega+\lambda'^2_{\ell}}{|\lambda'_{\ell}|}$, где $\lambda_{\ell}$ и $\lambda'_{\ell}$ – корни

многочленов Эйлера-Фробениуса $E_{2m-1}(z)$ и $E_{2m-3}(z)$ соответственно.

Причем, $\lambda_{m-1} < \lambda_m < ... < \lambda_1 < 0$ и $\lambda'_{m-2} < \lambda'_m < ... < \lambda'_1 < 0$ [3, С. 151]

Тогда можно записать, что

$$L(\omega) = |1-\cos\omega|^2 \frac{\prod_{\ell=1}^{m-1}\left(1-2\lambda_{\ell}\cos\omega+\lambda_{\ell}^2\right)}{\prod_{\ell=1}^{m-2}\left(1-2\lambda'_{\ell}\cos\omega+\lambda'^2_{\ell 1}\right)}.$$

Рассмотрим

$$L'(\omega) = \left((1-\cos\omega)^2\right)' \frac{\prod_{\ell=1}^{m-1}\left(1-2\lambda_{\ell}\cos\omega+\lambda_{\ell}^2\right)}{\prod_{\ell=1}^{m-2}\left(1-2\lambda'_{\ell}\cos\omega+\lambda'^2_{\ell}\right)} + (1-\cos\omega)\left(\frac{\prod_{\ell=1}^{m-1}\left(1-2\lambda_{\ell}\cos\omega+\lambda_{\ell}^2\right)}{\prod_{\ell=1}^{m-2}\left(1-2\lambda'_{\ell}\cos\omega+\lambda'^2_{\ell}\right)}\right)'$$

$$\tilde{L}'(\omega) = 2\sin\omega \frac{\prod_{\ell=1}^{m-1}\left(1-2\lambda_{\ell}\cos\omega+\lambda_{\ell}^2\right)}{\prod_{\ell=1}^{m-2}\left(1-2\lambda'_{\ell}\cos\omega+\lambda'^2_{\ell}\right)} + (1-\cos\omega)^2 \left(\frac{\prod_{\ell=1}^{m-1}\left(1-2\lambda_{\ell}\cos\omega+\lambda_{\ell}^2\right)}{\prod_{\ell=1}^{m-2}\left(1-2\lambda'_{\ell}\cos\omega+\lambda'^2_{\ell}\right)}\right)' = 0 \quad (3)$$

Очевидно, что $1-2\lambda_{\ell}\cos\omega+\lambda_{\ell}^2 > 0$ и $1-2\lambda'_{\ell}\cos\omega+\lambda'^2_{\ell} > 0$, а значит,

$\prod_{\ell=1}^{m-1}\left(1-2\lambda_{\ell}\cos\omega+\lambda_{\ell}^2\right) > 0$ и $\prod_{\ell=1}^{m-1}\left(1-2\lambda'_{\ell}\cos\omega+\lambda'^2_{\ell}\right) > 0$.

Таким образом, первое слагаемое в (3)

$2\sin\omega \frac{\prod_{\ell=1}^{m-1}\left(1-2\lambda_{\ell}\cos\omega+\lambda_{\ell}^2\right)}{\prod_{\ell=1}^{m-2}\left(1-2\lambda'_{\ell}\cos\omega+\lambda'^2_{\ell}\right)} = 0$, когда $\sin\omega = 0$, то есть когда $\omega = \pi$.

Второе слагаемое в (3) $(1-\cos\omega)^2 \left( \dfrac{\prod\limits_{\ell=1}^{m-1}\left(1-2\lambda_\ell\cos\omega+\lambda_\ell^{\ 2}\right)}{\prod\limits_{\ell=1}^{m-2}\left(1-2\lambda'_\ell\cos\omega+\lambda'_\ell{}^{2}\right)} \right)'$ обращается в ноль, когда, либо $(1-\cos\omega)^2=0$, ( т.е. $\omega=\pi$ ), либо $M(\omega)=\left( \dfrac{\prod\limits_{\ell=1}^{m-1}\left(1-2\lambda_\ell\cos\omega+\lambda_\ell^{\ 2}\right)}{\prod\limits_{\ell=1}^{m-2}\left(1-2\lambda'_\ell\cos\omega+\lambda'_\ell{}^{2}\right)} \right)'=0$

Так как $1-2\lambda_\ell\cos\omega+\lambda_\ell^{\ 2}>0$ и $1-2\lambda'_\ell\cos\omega+\lambda'_\ell{}^{2}>0$, то достаточно рассмотреть числитель производной функции $M(\omega)$.

$\left(1-2\lambda_\ell\cos\omega+\lambda_\ell^{\ 2}\right)'\left(1-2\lambda'_\ell\cos\omega+\lambda'_\ell{}^{2}\right)-\left(1-2\lambda'_\ell\cos\omega+\lambda'_\ell{}^{2}\right)'\left(1-2\lambda_1\cos\omega+\lambda_1^{\ 2}\right)=$
$2\lambda_\ell\sin\omega\left(1-2\lambda'_\ell\cos\omega+\lambda'_\ell{}^{2}\right)-\left(2\lambda'_\ell\sin\omega\right)\left(1-2\lambda_\ell\cos\omega+\lambda_\ell^{\ 2}\right)=$
$2\sin\omega\left(\lambda'_\ell-\lambda_{\phantom{\ell}}\right)\left(1-\lambda'_\ell\lambda_\ell\right)=0$

Таким образом, второе слагаемое в (3) обращается в $0$, только когда $\omega=\pi$. Рассмотрим $\left(\lambda'_\ell-\lambda_{\phantom{\ell}}\right)\left(1-\lambda'_\ell\lambda_\ell\right)$. Если это произведение больше нуля, то можно утверждать, что знак второго слагаемого в (3) в окрестности точки $\omega=\pi$ зависит от знака $\sin\omega$.

И $\max\limits_{\omega} L(\omega)$ будет достигаться в точке $\omega_0=\pi$.

Очевидно, что корни полиномов Эйлера–Фробениуса $|\lambda'_1|>|\lambda_{1-1}|$.

Значит,

$$\left(\lambda'_\ell-\lambda_{\phantom{\ell}}\right)\left(1-\lambda'_\ell\lambda_\ell\right)>\left(\lambda_{\ell-1}-\lambda_1\right)\left(1-\lambda'_\ell\lambda_\ell\right)>\left(\lambda_{\ell-1}-\lambda_\ell\right)\left(1-\lambda_{\ell-1}\lambda_\ell\right)>0$$

Таким образом,

$L'(\omega)=0$ в точке $\omega_0=\pi$ и в этой точке функция $L(\omega)$ достигает своего максимума.

То есть,

$$\max\limits_{\omega} L(\omega)=2|1-\cos\pi|\dfrac{\sum\limits_\ell\left|\hat{N}_{m-1}(\pi+2\pi\ell)\right|^2}{\sum\limits_\ell\left|\hat{N}_m(\pi+2\pi\ell)\right|^2}=2\dfrac{\sum\limits_\ell\left|\hat{N}_{m-1}(\pi+2\pi\ell)\right|^2}{\sum\limits_\ell\left|\hat{N}_m(\pi+2\pi\ell)\right|^2}.$$

Функцию $\sum\limits_\ell\left|\hat{N}_m(\omega+2\pi\ell)\right|^2$ в общем виде можно записать

$$\sum_{\ell}\left|\hat{N}_m(\omega+2\pi\ell)\right|^{2m} = \frac{2^{2m+2}\sin^{2m+2}\frac{\omega}{2}}{\sum_{\ell}|\omega+2\pi\ell|^{2m+2}}$$

соответственно

$$\sum_{\ell}\left|\hat{N}_{m-1}(\omega+2\pi\ell)\right|^{2m-2} = \frac{2^{2m}\sin^{2m}\frac{\omega}{2}}{\sum_{\ell}|\omega+2\pi\ell|^{2m}}.$$

Тогда,

$$\max_{\omega} L(\omega) = 2\frac{\sum_{\ell}\left|\hat{N}_{m-1}(\pi+2\pi\ell)\right|^{2m}}{\sum_{\ell}\left|\hat{N}_m(\pi+2\pi\ell)\right|^{2m+2}} = 2\frac{2^{2m}\sin^{2m}\frac{\pi}{2}}{\sum_{\ell}|\pi+2\pi\ell|^{2m}} \frac{\sum_{\ell}|\pi+2\pi\ell|^{2m+2}}{2^{2m+2}\sin^{2m+2}\frac{\pi}{2}}$$

$$= \frac{1}{\sum_{\ell}|\pi+2\pi\ell|^{2m}} \frac{\sum_{\ell}|\pi+2\pi\ell|^{2m+2}}{2\sin^2\frac{\pi}{2}} = \pi^2\frac{\sum_{\ell}|\pi+2\pi\ell|^{2m+2}}{\sum_{\ell}|\pi+2\pi\ell|^{2m}} = \pi^2\frac{\sum_{\ell}|1+2\pi\ell|^{2m+2}}{\sum_{\ell}|1+2\pi\ell|^{2m}} \quad (4)$$

Величину в правой части (4) удобно записать с помощью константы Фавара:

$$K_m = \frac{4}{\pi}\sum_{\ell=0}^{\infty}\frac{(-1)^{\ell(m+1)}}{(1+2\ell)^{m+1}}, \quad m=0,1,2\ldots \quad [2, 64\text{-}65]$$

Таким образом (4) можно переписать: $\max_{\omega} L(\omega) = \frac{\pi^2 K_{2m-1}}{K_{2m+1}}$.

Тогда, согласно неравенству (2), имеем

$$\|s'\|_2 \leq \pi\sqrt{\frac{K_{2m-1}}{K_{2m+1}}}\|s\|_2.$$

Если же рассмотреть целочисленные сдвиги сплайнов, то есть сплайны, вида $s(hx+k)$, ( где $h$ - шаг, $h \in N$ ) тогда можно записать, что

$$\|s'\|_2 \leq \pi h\sqrt{\frac{K_{2m-1}}{K_{2m+1}}}\|s\|_2.$$

Теперь, воспользовавшись методом математической индукции, получим оценку нормы $k$-ой производной для непериодического сплайна в пространстве $L_2$:

$\forall k \in Z, \ k \leq m, \ \forall h \in N$ справедливо

$$\|s^{(k)}\|_2 \leq (\pi h)^k\sqrt{\frac{K_{2(m-k)+1}}{K_{2m+1}}}\|s\|_2.$$

Константа в правой части полученного неравенства неулучшаема.

Чтобы в этом убедиться, достаточно рассмотреть неулучшаемость константы в правой части неравенства (2).

Положим
$$\left|m_s(\omega)\right|^2 = \Phi_n(\omega - \omega_0),$$

где $\Phi_n(\omega)$ - ядро Фейера, порядка $n$, $\omega_0$ - точка, которая реализует максимум в правой части неравенства (2). Заметим, что $\dfrac{1}{2\pi}\int_0^{2\pi}\Phi_n(\omega)d\omega = 1$.

Положим далее $\hat{s}(x) = m_s(\omega)\hat{N}_m(\omega)$. Заметим также, что $\forall \ell \in Z$, $m = 0,1,2\ldots$

$$\|s(x)\|_2^2 = \int_0^{2\pi}\left|m_s(\omega)\right|^2 \sum_\ell \left|\hat{N}_m(\omega+2\pi\ell)\right|^2 d\omega.$$

Таким образом, будем иметь

$$\|s'(x)\|_2^2 = \frac{1}{2\pi}\int_R \left|\sum_\gamma c_\gamma e^{-i\gamma x}\hat{N}_{m-1}(\omega)(1-e^{-i\omega})\right|^2 d\omega =$$

$$\frac{1}{2\pi}\int_0^{2\pi}\sum_\ell \left|m_s(\omega)\right|^2 \left|\hat{N}_{m-1}(\omega+2\pi\ell)(1-e^{-i\omega})\right|^2 d\omega =$$

$$\frac{1}{2\pi}\int_0^{2\pi} \frac{\sum_\ell \left|i\omega+2\pi\ell\right|^2 \left|\hat{N}_m(\omega+2\pi\ell)\right|^2}{\sum_\ell \left|\hat{N}_m(\omega+2\pi\ell)\right|^2}\left|\Phi_n(\omega-\omega_0)\right|\sum_\ell \left|\hat{N}_{m-1}(\omega+2\pi\ell)\right|^2\left|1-e^{-i\omega}\right|^2 d\omega$$

$$\xrightarrow[n\to\infty]{} \frac{\left|\hat{N}_{m-1}(\omega+2\pi\ell)\right|^2\left|1-e^{-i\omega}\right|^2}{\sum_\ell \left|\hat{N}_m(\omega+2\pi\ell)\right|^2}\|s(\omega)\|_2^2.$$

Таким образом, $\dfrac{\|s'(\omega)\|_2^2}{\|s(\omega)\|_2^2} \to \max_\omega \dfrac{\left|\hat{N}_{m-1}(\omega+2\pi\ell)\right|^2\left|1-e^{-i\omega}\right|^2}{\sum_\ell \left|\hat{N}_m(\omega+2\pi\ell)\right|^2}$ .